\documentclass[11pt, dvips, twoside]{amsart} 

\usepackage{amssymb}
\usepackage{psfrag,epsfig}

\theoremstyle{plain}
\newtheorem{theorem}{Theorem}[section]
\newtheorem{lemma}[theorem]{Lemma}
\newtheorem{proposition}[theorem]{Proposition}
\newtheorem{corollary}[theorem]{Corollary}
\newtheorem{remark}[theorem]{Remark}
\newtheorem{remarks}[theorem]{Remarks}

\newtheorem{examples}[theorem]{Examples}
\newtheorem{definition}[theorem]{Definition}

\newtheorem{problem}[theorem]{Open Problem}
\newtheorem{problems}[theorem]{Open Problems}

\newcommand{\proofend}{\hspace*{\fill} $\Box$\\}
\newcommand{\diam}{\hspace*{\fill} $\Diamond$\\}

\def\s{\smallskip}
\def\m{\medskip}
\def\eps{\varepsilon}

\def\Vol{\operatorname {vol}\:\!}

\def\Diffc0{\operatorname{Diff^c_0}}

\def\Sympc0{\operatorname{Symp^c_0}}

\def\Int{\operatorname{Int}}

\def\length{\operatorname{length}}
\def\supp{\operatorname{supp}}

\def\Crit{\operatorname{Crit}}

\def\Diam{\operatorname{Diam}}
\def\G{\operatorname{G}}
\def\HZ{\operatorname{HZ}}
\def\Oh{\operatorname{Oh}}
\def\PSS{\operatorname{PSS}}
\def\V{\operatorname{V}}
\def\H{\operatorname{H}}
\def\HM{\operatorname{HM}}
\def\HF{\operatorname{HF}}
\def\CF{\operatorname{CF}}
\def\CM{\operatorname{CM}}

\def\ga{\alpha}
\def\gb{\beta}
\def\gg{\gamma}
\def\gd{\delta}
\def\eps{\epsilon}

\def\gf{\varphi}

\def\gl{\lambda}
\def\go{\omega}
\def\gs{\sigma}

\def\ca{{\mathcal A}}

\def\cd{{\mathcal D}}

\def\cf{{\mathcal F}}

\def\ch{{\mathcal H}}

\def\cl{{\mathcal L}}
\def\cm{{\mathcal M}}

\def\cp{{\mathcal P}}
\def\cs{{\mathcal S}}

\def\CC{\mathbb{C}}

\def\RR{\mathbb{R}}
\def\ZZ{\mathbb{Z}}

\def\pp{\partial}

\def\ra{\rightarrow}

\def\ni{\noindent}
\def\b{\bigskip}
\def\m{\medskip}

\def\id{\mbox{id}}

\def\proof{\noindent {\it Proof. \;}}

\begin{document}

\title{Energy capacity inequalities via an action selector}
$   $ \\
$   $ \\

\author{Urs Frauenfelder}
\address{(U.\ Frauenfelder) Department of Mathematics, Hokkaido University,
Sapporo 060-0810, Japan}
\email{urs@math.sci.hokudai.ac.jp}
\author{Viktor Ginzburg}
\address{(V.\ Ginzburg) Department of Mathematics, UC Santa
Cruz, CA 95064, USA}
\email{ginzburg@math.ucsc.edu}
\author{Felix Schlenk}
\address{(F.\ Schlenk) Mathematisches Institut,
Universit\"at Leipzig, 04109 Leipzig, Germany}
\email{schlenk@math.uni-leipzig.de}

\date{\today}
\thanks{This work is supported in part by the Swiss National Science
Foundation (UF), the US National Science Foundation
and the UCSC faculty
research funds (VG), and the Deutsche Forschungsgemeinschaft (FS)}

\begin{abstract}
An action selector for a symplectic manifold $(M,\go)$
associates with each compactly supported Hamiltonian function
$H$ on $M$ an action value of $H$ in a suitable way.
Action selectors are known to exist for a broad class of 
symplectic manifolds.
We show how the existence of an 
action selector leads to sharp energy capacity
inequalities between the Gromov width, the Hofer--Zehnder
capacity, and the displacement energy. We also obtain sharp
lower bounds for the smallest action of a closed characteristic
on contact type hypersurfaces.
\end{abstract}

\maketitle

\markboth{
{\rm Energy capacity inequalities via an action selector}
}{{}} 

\tableofcontents

\section{Introduction and main results}

\ni
Consider an arbitrary symplectic manifold $(M,\go)$.
We set $I = [0,1]$ and denote, for each subset $A$ of $M$,
by $\ch (I \times A)$ the set of smooth functions $H \colon I
\times M \ra \RR$ whose support is compact and contained in
$I \times \Int A$.
We abbreviate $\ch = \ch (I \times M)$, and denote by $\ch (A)$
the set of those functions in $\ch (I \times A)$ which do not
depend on $t \in I$.
The Hamiltonian vector field of $H \in \ch$ defined by 
\begin{equation}  \label{def:H}
\go \left( X_{H_t}, \cdot \right) \,=\, -d H_t \left( \cdot \right)
\end{equation}
generates a flow $\gf_H^t$ with time-$1$-map $\gf_H$.
For $H \in \ch$, the set of contractible $1$-periodic orbits of
$\gf_H^t$ is denoted $\cp^\circ (H)$.
Given $x \in \cp^\circ (H)$, let $\cd (x)$ be the set of smooth
discs $\bar{x} \colon D^2 = \left\{ z \in \CC \mid |z| \le 1 \right\}
\ra M$ satisfying $\bar{x}(e^{it})=x(t)$.
We shall identify the map $\bar{x}$ with its oriented image and
write
$\int_{\bar x} \go = \int_{D^2} \left( \bar x \right)^* \go$.
For $H \in \ch$, the action functional $\ca_H$ on 
$\bar{\cp}^\circ (H) = \left\{ (x,\bar{x}) \mid x \in \cp^\circ
(H), \bar{x} \in \cd (x) \right\}$ is defined as
\[
\ca_H \left( \bar{x} \right) \,=\, 
      -\int_{\bar x} \go + \int_0^1 H \left(t, x(t) \right) dt ,
\]
and its {\it action spectrum}\, is
\[
\Sigma^\circ (H) \,=\, \left\{ \ca_H \left( \bar{x} \right)
\mid (x,\bar{x}) \in \bar\cp^\circ (H) \right\} .
\]
Note that for $a\in \Sigma^\circ (H)$ the set
$a+\omega(\pi_2(M))$ also belongs to $\Sigma^\circ (H)$. As a 
consequence, $\Sigma^\circ (H)$ need not be closed and is dense
in $\RR$ if $\omega(\pi_2(M))$ is dense.
For $H \in \ch$ we abbreviate 
\[
E^+(H) \,=\, \int_0^1 \max_{x \in M} H(t,x) \,dt ,
\]
and we recall that
for $H,K \in \ch$ the composition $\gf_H \circ \gf_K$ 
is generated by
\[
\left( H \# K \right)(t,x) \,=\, 
H(t,x) + K \left( t,\left( \gf^t_H \right)^{-1}(x) \right) .
\]
We say that a Hamiltonian $H \in \ch (M)$ is {\it simple}\, and
write $H \in \cs (M)$ if
\begin{itemize}
\item[(P1)] 
  $H \ge 0$,
\s
\item[(P2)] 
  $H |_U = \max H$ for some open non-empty set $U \subset M$,
\s
\item[(P3)] 
  the only critical values of $H$ are $0$ and $\max H$.
\end{itemize}
We emphasize that simple Hamiltonians are, as is clear from (P3),
``normalized'' to have minimum equal to $0$.
Furthermore, we say that $H \in \cs (M)$ is {\it Hofer--Zehnder
admissible}\, and write $H \in \cs_{\HZ}^\circ(M)$ if 
the flow $\gf_H^t$ has no non-constant, contractible in $M$, $T$-periodic 
orbit with period $T \le 1$.

\subsection{Axioms for an action selector}
A {\it weak action selector}\, $\gs$ for $(M,\go)$ is a map
$\gs \colon \ch \ra \RR$ satisfying
the following axioms.
{\rm
\begin{itemize}
\item[(AS1)] 
  $\gs(H) \in \Sigma^\circ(H)$\, for all $H \in \ch$;
\s
\item[(AS2)] 
  $\gs(H) >0$ for all $H \in \cs (M)$ with $H \not\equiv 0$;
\s
\item[(AS3)]
  $\gs(H) \le E^+(H)$\, for all $H \in \ch$;
 \s
\item[(AS4)]
 $\gs$ is continuous with respect to the $C^0$-topology on $\ch$;
 \s
\item[(AS5)] 
 $\gs \left(H \# K \right) \le \gs (H) + E^+(K)$\, for all $H,K
 \in \ch$.
\end{itemize}
An {\it action selector}\, $\gs$ for $(M,\go)$ is a weak action
selector which {\it in addition} to (AS2) satisfies
\begin{itemize}
\item[(AS2$^+$)] 
  $\gs(H) = \max H$\, for all $H \in \cs_{\HZ}^\circ (M)$.
\end{itemize}

A symplectic manifold $(M, \go)$ is called {\it weakly exact}\,
if $[\go]$ vanishes on $\pi_2(M)$.
\begin{remark}  \label{r:weak}
{\rm
If $(M, \go)$ is weakly exact, a weak action selector for
$(M, \go)$ is an action selector for $(M, \go)$.
}
\end{remark}
\proof
Given $H \in \cs_{\HZ}^\circ (M)$ on a weakly exact
symplectic manifold $(M, \go)$, we have $\Sigma^\circ (H) = \{
0, \max H \}$, and hence axioms (AS1) and (AS2) imply $\gs (H) =
\max H$.
\proofend
\begin{remark}  \label{r:AS3} 
{\rm
If $\gs$ is an action selector for $(M, \go)$, 
then (AS3) is a consequence of the other axioms.
}
\end{remark}
\proof
Axioms (AS5) and (AS2$^+$) applied to $H=0$ yield 
\[
\gs (K) \,=\, \gs \left( 0 \# K \right) \,\le\, \gs (0) + E^+(K) \,=\,
E^+(K)
\] 
for all $K \in \ch$.
\proofend

An {\it exhaustion}\, of a symplectic manifold 
$(M,\go)$ is an increasing sequence of submanifolds 
$M_i \subset M$ exhausting $M$, that is,
\[
M_1 \subset M_2 \subset \dots \subset M_i \subset \dots \subset M
\qquad \text{and} \qquad
\bigcup_i M_i = M . 
\]  

By a (weak) action selector for an exhaustion $(M_i)$ of
$(M,\go)$ we mean a collection $\gs_i$ of (weak) action
selectors for $(M_i, \go)$.
While a (weak) action selector for $(M,\go)$ obviously restricts to a (weak)
action selector for any of its exhaustions, the
restriction of $\gs_{i+1}$ to $M_i$ is not assumed to agree
with $\gs_i$, whence it is unclear whether
every (weak) action selector for an exhaustion of $(M,\go)$ fits
together to a (weak) action selector for $(M,\go)$;
see 3.\ (ii) and (iii) in Appendix~\ref{app:0}.
The existence of a (weak) action selector for exhaustions of
$(M,\go)$ as defined here is sufficient for our purposes in this paper and
can sometimes be established even when it is unknown whether a global
(weak) action selector for $(M,\go)$ exists; see
3.\ (iii) in Appendix~\ref{app:0}.

\b
\ni
{\bf Examples of (weak) action selectors.}

\m
\ni
{\bf 1. $\left( \RR^{2n}, \go_0 \right)$.}
Let $\go_0 = \sum_i dq_i \wedge dp_i$ be the standard symplectic
form on $\RR^{2n}$. Action selectors $\gs_{\V}$ and $\gs_{\HZ}$ 
for $( \RR^{2n}, \go_0)$ have been constructed by Viterbo in
\cite{V} and by Hofer and Zehnder in \cite{H, HZ}.
The constructions of $\gs_{\V}$ and $\gs_{\HZ}$ are outlined in 
Appendix~\ref{app:0}.

\m
\ni
{\bf 2. Weakly exact closed symplectic manifolds.}
After work on selector-like invariants for standard cotangent bundles 
over a closed base by Oh \cite{Oh1,Oh2},
an action selector $\gs_{\PSS}$ was constructed
for weakly exact closed symplectic manifolds by Schwarz in \cite{Sch} by
making use of the Piunikhin--Salamon--Schwarz isomorphism.
Examples of such manifolds are (products of) closed surfaces of
positive genus. 
The construction of $\gs_{\PSS}$ is outlined in Appendix~\ref{app:0}.

\m
\ni
{\bf 3. Weakly exact convex symplectic manifolds.}
A compact symplectic manifold $(M, \go)$ with
boundary $\pp M$ is said to be {\it convex}\, if 
there exists a Liouville vector field $X$ 
(i.e., $\cl_X \go = d \iota_X \go \stackrel{!}{=} \go$) 
which is defined near $\pp M$ and is everywhere transverse to $\pp M$, 
pointing outward.
A non-compact symplectic manifold $(M, \go)$ is {\it convex}\,
if it admits an exhaustion by compact convex submanifolds.
Examples of weakly exact convex symplectic manifolds are
cotangent bundles $\left( T^*B, \go_0 \right)$
over a closed base $B$ endowed with the standard symplectic form
$\go_0 = \sum_i dp_i \wedge dq_i$, and unit-ball bundles
therein, and, more generally, Stein manifolds and Stein domains,
see 3.\ (i) in Appendix~\ref{app:0}.
Other examples are twisted cotangent bundles over closed
orientable surfaces of genus at least $2$, see \cite{FS,Gi0}. 
Building on \cite{Sch} and \cite{V2}, an action selector
$\gs_{\PSS}$ for (exhaustions of) weakly exact convex symplectic
manifolds was constructed in \cite{FS}.

\m
\ni
{\bf 4. Rational strongly semi-positive closed symplectic manifolds.}
A $2n$-di\-men\-sion\-al symplectic manifold $(M,\go)$ is {\it
strongly semi-positive}\, 
if it satisfies one of the following conditions.
{\rm
\begin{itemize}
\item[(SP1)] 
  $\go (A) = \gl\, c_1 (A)$ for every $A \in \pi_2(M)$ where $\gl \ge 0$;
\s
\item[(SP2)] 
  $c_1(A)=0$ for every $A \in \pi_2(M)$;
\s
\item[(SP3)] 
  the minimal Chern number $N \ge 0$ defined by 
  $c_1 \left( \pi_2(M) \right) = N \ZZ$ is at least $n-1$.
\end{itemize}
Here, $c_1 = c_1 (\go)$ is the first Chern class of any almost
complex structure on $TM$ compatible with $\go$, see
\cite[Section~4.1]{MS}.
Examples are symplectic manifolds of dimension $\le 4$, complex
Grassmannians, and Calabi--Yau manifolds.
A symplectic manifold $(M, \go)$ is {\it rational}\, if $\go
\left( \pi_2(M)\right)$ is a discrete subset of $\RR$.
The action selector $\gs_{\PSS}$ from \cite{Sch} for weakly
exact closed symplectic manifolds has been extended to 
a {\it weak}\, action selector $\gs_{\PSS}$ for rational
strongly semi-positive
closed symplectic manifolds in \cite{MS3}.

\m
\ni
{\bf 5. Rational closed symplectic manifolds.}
A construction of an action selector $\gs_{\Oh}$ for all
rational closed symplectic
manifolds was proposed by Oh in \cite{Oh3,Oh4,Oh5}.
\diam
\begin{remarks}\  
{\rm
{\bf 1.}
In view of the action selector $\gs_{\Oh}$, our
concept of a weak action selector appears to be superfluous.
The verification of (AS2$^+$) for $\gs_{\Oh}$ is, however,
difficult already for rational strongly semi-positive closed symplectic
manifolds. We thus find it interesting to see which results can
be formally derived from the existence of a weak action selector.

\s
\ni
{\bf 2.}
We have chosen a minimal set of axioms for a (weak) action selector
required for its applications that we have in mind: The existence of
a (weak) action selector for a given symplectic manifold will imply
Theorems~1 and 2 below. 
The selector-like invariants constructed in 
\cite{Oh1,Oh2} and the (weak)
action selectors $\gs_{\V}$, $\gs_{\HZ}$, $\gs_{\PSS}$ and
$\gs_{\Oh}$ have many further 
properties, and they lead to many other results in Hamiltonian dynamics 
\cite{FS,FS2,Gi,HZ,Oh1,Oh2,Oh4,Oh5,P1,Sch,V} 
as well as to some insight into the algebraic structure of the groups
of Hamiltonian and symplectic diffeomorphisms of certain closed
symplectic manifolds \cite{E,EP,P1}.
An example of an additional property 
is sub-additivity:
\begin{equation*}  
\text{$\sigma(H \# K)\le \sigma(H)+\sigma(K)$\, for all $H,K \in
\ch$}.
\end{equation*}
This together with (AS3) is stronger than the axioms (AS3) and 
(AS5), which are sufficient for our purposes.
Sub-additivity holds for the selectors 
$\gs_{\V}$, $\gs_{\PSS}$ and $\gs_{\Oh}$, but has not been established for
$\gs_{\HZ}$.
An additional property shared by all known (weak)
action selectors and useful for intuition is the monotonicity property
\begin{equation*}  
\text{$\gs(H) \le \gs (K)$\, for all $H,K \in \ch$ with $H \le K$.}
\end{equation*}
}
\end{remarks}

\begin{problems}\
{\rm
{\bf 1.}
It is not known whether every (exhaustion of a non-compact)
symplectic manifold admits a (weak) action selector.
The simplest symplectic manifold for which no weak action selector
has yet been constructed is 
the convex strongly semi-positive symplectic manifold 
$\left( T^* S^2, \go_0 + \pi^* \gs \right)$, where $\gs$ is an
area form on $S^2$.

\s
\ni
{\bf 2.}
Is every weak action selector an action selector?

\s
\ni
{\bf 3.}
Do axioms (AS1)--(AS5) uniquely determine (weak) action selectors?
This would imply that every (weak) action selector for an
exhaustion of $(M,\go)$ fits together to a (weak) action
selector for $(M,\go)$.
More specifically, the action selectors $\gs_{\V}$,
$\gs_{\HZ}$, $\gs_{\PSS}$ and $\gs_{\Oh}$ are all defined by a
variational procedure (see Appendix~\ref{app:0}), 
and so it should be possible to compare them:
Do $\gs_{\V}$, $\gs_{\HZ}$ and $\gs_{\PSS}$ agree on $\left(\RR^{2n},
\go_0 \right)$?
Does $\gs_{\PSS}$ agree with $\gs_{\Oh}$ on weakly exact or 
strongly semi-positive closed symplectic manifolds?
}         
\end{problems}

\subsection{A sharp energy capacity inequality}

Given a symplectic manifold $(M,\go)$ one can associate to each
subset $A$ of $M$ various symplectic invariants.

\s
\ni
{\bf 1. The Gromov width.} 
The Gromov width of $A$ is defined
as
\[
c_{\G} (A) \,=\, \sup \left\{ \pi r^2 \mid B^{2n}(r) \text{
symplectically embeds into } (M, \go) \right\} .
\]
Here, $B^{2n}(r)$ denotes the open ball in $\left( \RR^{2n},
\go_0 \right)$ of radius $r$.
The Gromov width, which was introduced by Gromov in \cite{Gr}, 
measures the symplectic size of $(M, \go)$ in a
geometric way; it corresponds to the injectivity radius of a
Riemannian manifold. Darboux's theorem states that every point 
of a symplectic manifold has an open neighborhood which is
symplectomorphic to some ball $B^{2n}(r)$, and hence $c_{\G} (A)$
vanishes only if $A$ has empty interior.

\s
\ni
{\bf 2. Hofer--Zehnder capacities.} 
Hofer--Zehnder capacities, in contrast with the Gromov width, measure the 
symplectic size of a set from
the perspective of Hamiltonian dynamics on this set.
We consider two variants.
For each subset $A \subset M$ let $\cs (A)$ be the
set of simple functions in $\ch (A)$.
We say that a function $H \in \cs (A)$ is {\it $\HZ$-admissible}\, if the flow
$\gf_H^t$ has no non-constant $T$-periodic orbit with period $T \le 1$,
and as before $H \in \cs (A)$ is {\it $\HZ^\circ$-admissible}\, 
if the flow $\gf_H^t$ has no non-constant $T$-periodic orbit 
with period $T \le 1$ which is contractible in $M$.
Set
\begin{eqnarray*}
\cs_{\HZ} (A) &=& \left\{ H \in \cs (A) \mid H \text{ is
HZ-admissible} \right\}, \\[0.2em]
\cs_{\HZ}^\circ (A,M) &=& \left\{ H \in \cs (A) \mid H \text{ is
$\HZ^\circ$-admissible} \right\}.
\end{eqnarray*}
Following \cite{HZ2,HZ} and \cite{Lu,Sch}, we define 
the Hofer--Zehnder capacity and the $\pi_1$-sensitive Hofer--Zehnder
capacity of $A \subset (M, \go)$ as
\begin{eqnarray*}
c_{\HZ} (A) &=& \sup \left\{ \max H \mid H \in \cs_{\HZ}
(A) \right\}, \\[0.2em]
c_{\HZ}^\circ (A,M) &=& \sup \left\{ \max H \mid H \in
\cs_{\HZ}^\circ (A,M) \right\} .
\end{eqnarray*}
Both $c_{\HZ} (A)$ and $c_{\HZ}^\circ (A,M)$ vanish if
and only if $A$ has empty interior.
\begin{remarks}\  \label{r:HZ}
{\rm
{\bf 1.}
The definition of the Hofer--Zehnder capacities $C_{\HZ} (A)$ and
$C_{\HZ}^\circ (A,M)$ 
given in \cite{HZ2,HZ} and \cite{Lu,Sch} uses the larger
class 
\[
\cf (A) \,=\, \left\{ H \in \ch \left( A \right) \mid
              H \text{ satisfies (P1) and (P2)} \right\} .
\]
Of course, $c_{\HZ} (A) \le C_{\HZ} (A)$ 
and $c_{\HZ}^\circ (A,M) \le C_{\HZ}^\circ (A,M)$. 
It follows from \cite{HZ2} that equalities hold for convex subsets
of $\left( \RR^{2n}, \go_0 \right)$, and we do not know examples with 
$c_{\HZ} (A) < C_{\HZ} (A)$ or $c_{\HZ}^\circ (A,M) < C_{\HZ}^\circ (A,M)$. 

\s
\ni
{\bf 2.}
The main virtue of the Hofer--Zehnder capacity $C_{\HZ} (A)$ is that
$C_{\HZ} (A) < \infty$ implies almost existence of periodic orbits near
any compact regular energy level of an autonomous Hamiltonian system on
$A$, and similarly for $C_{\HZ}^\circ (A,M)$. 
As we shall show in Appendix~\ref{app:1}, 
this continues to hold for $c_{\HZ} (A)$ and $c_{\HZ}^\circ (A,M)$.

\s
\ni
{\bf 3.}
Corollary~1 below remains valid for
$C_{\HZ}$ and $C_{\HZ}^\circ$ when (AS2$^+$) is required to hold
for the larger classes $\cf_{\HZ}(A)$ and $\cf_{\HZ}^\circ(A,M)$
of Hofer--Zehnder admissible functions from $\cf (A)$. The
action selectors $\gs_{\V}$, $\gs_{\HZ}$, $\gs_{\PSS}$ and $\gs_{\Oh}$ from
Examples~$1$, $2$, $3$ and $5$ do satisfy this stronger axiom.
}
\end{remarks}

\s
\ni
{\bf 3. Displacement energy.} 
An invariant with both geometric and dynamical features is the
{\it displacement energy}\, introduced in \cite{Hofer90,LM1}.
The Hofer norm $\| H \|$ of $H \in \ch$ is defined as
\[
\| H \| \,=\, \int_0^1 
\left( \sup_{x \in M} H(t,x) - \inf_{x \in M} H(t,x) \right) dt ,
\]
and the displacement energy $e(A,M) = e(A,M,\go) \in [0,\infty]$ 
is defined as
\[
e(A,M) \,=\, \inf \left\{ \left\| H \right\| \mid H \in
\ch,\,
\gf_H (A) \cap A = \emptyset \right\} 
\]
if $A$ is compact and as
\[
e(A,M) \,=\, \sup \left\{ e(K,M) \mid K \subset A \text{ is compact} \right\} 
\]
for a general subset $A$ of $M$.

\m
Since the invariants $c_{\G}$, $c_{\HZ}$ and $c_{\HZ}^\circ$, and
$e$
are defined in different ways, relations
between them lead to many applications.
It is easy to see that $c_{\G} (A) \le c_{\HZ} (A)$, see e.g.\
\cite{HZ}, and it follows from definitions that
$c_{\HZ} (A) \le c_{\HZ}^\circ (A,M)$. In order to compare
$c_{\HZ}^\circ (A,M)$ with $e(A,M)$, we introduce further
invariants. Following \cite{V}, we make the
\begin{definition}
{\rm
Assume that $(M, \go)$ admits a weak action selector $\gs$. For each
subset $A$ of $M$ the {\it spectral capacities}\, $c_\gs (A,M)$ and $c^\gs (A,M)$
are defined as
\begin{eqnarray*}
c_\gs (A,M) &=& \sup \left\{ \gs (H) \mid H \in \cs (A) \right\} ,   \\
c^\gs (A,M) &=& \sup \left\{ \gs (H) \mid H \in \ch(I \times A) \right\} .
\notag
\end{eqnarray*}
}
\end{definition}
\ni
Of course, $c_\gs (A,M) \le c^\gs (A,M)$, 
and $c_\gs (A,M) = c^\gs (A,M)$ provided that $\gs$ is monotone and that $A
\subsetneq M$ if $M$ is closed.
If $M$ is closed, $c^\gs(M,M) = \infty$, because then
$\ch(I\times M)$ contains the constant functions and
$\sigma(\mathit{const})=\mathit{const}$. 
(Indeed, $\sigma (0) \le 0$ by (AS3) and $\gs (0) \ge 0$ by
(AS2) and (AS4), so that 
$\gs (0)=0$. This, together with (AS1), (AS4) and the fact that $\go \left(
\pi_2(M) \right)$ is countable, implies that 
$\sigma(\mathit{const})=\mathit{const}$.) 
Since simple Hamiltonians have minimum $0$, this argument does not apply to
$c_\gs$. However, for all closed $M$ for which we know $c_\gs(M,M)$,
even this capacity is infinite.

\b
The following result, proved in Section \ref{sec:proofs}, is an elaboration 
of an observation in \cite{Gi3}.
\\
\\
{\bf Theorem~1.}
{\it Consider a symplectic manifold $(M,\go)$ and an arbitrary subset $A$ of $M$.

\s
\ni
(i) If $\gs$ is a weak action selector for $(M,\go)$, then
\[
c_\gs (A,M) \,\le\, c^\gs (A,M) \,\le\, e (A,M) .
\]

\ni
(ii) If $\gs$ is an action selector, 
then $c_{\HZ}^\circ (A,M) \le c_\gs (A,M)$, so that
\[
c_{\G}(A) \,\le\, c_{\HZ} (A) \,\le\, c_{\HZ}^\circ (A,M) \,\le\, 
c_\gs (A,M) \,\le\, c^\gs (A,M) \,\le\, e (A,M) .
\]
}
\\
{\bf Corollary~1.}
{\it
Assume that an exhaustion of $(M,\go)$ admits an action selector. 
Then
\[
c_{\G}(A) \,\le\, c_{\HZ} (A) \,\le\, c_{\HZ}^\circ (A,M) \,\le\, 
e (A,M) 
\]
for every subset $A$ of $M$.
}

\b
\noindent{\it Proof of Corollary~1.}
The first two inequalities in the corollary are clear.
Let $(M_i)$ be an exhaustion of $(M,\go)$ with action selectors
$\gs_i$.
For every $i$ with $A \subset M_i$, Theorem~1\,(ii) yields
\[
c_{\HZ}^\circ \left( A,M_i \right) \,\le\, c_\gs \left( A,M_i
\right) \,\le\, e \left( A,M_i \right) ,
\]
and so we readily find 
$c_{\HZ}^\circ \left( A,M \right) \le e \left( A,M \right)$.
\proofend
\begin{remark}\  \label{r:sharp}
{\rm
{\bf 1.}
Theorem~1 and Corollary~1 are sharp. 
Indeed, let $\gf \colon B^{2n}(3r) \hookrightarrow (M,\go)$ 
be a Darboux ball. For $A=\gf \left( B^{2n}(r) \right)$ we have
$c_G(A) = e (A,M) = \pi r^2$.

\s
\ni
{\bf 2.}
It would be interesting to know whether in the situation of
Theorem~1\,(ii) there exists an example
with $c_{\HZ}^\circ (A,M) < c_\gs (A,M)$.
A more specific version of this problem is posed in \ref{open:<}
below.

\s
\ni
{\bf 3.}
(i)
The assertion of Corollary~1 was first obtained by Hofer,
\cite{H}, for $(\RR^{2n}, \go_0 )$, 
see also \cite[Section~5.5]{HZ};
in fact, our axioms for a (weak) action selector are extracted
from \cite{HZ}, and our proof of Theorem~1 closely follows
\cite{HZ}.
Later on, the energy-capacity inequality  
$c_{\HZ}^\circ \left( A, M \right) \le 2\,e \left( A, M \right)$
was established for subsets of weakly exact 
symplectic manifolds which are closed 
\cite{Sch} or convex \cite{FS},
and it was pointed out in \cite{Gi3}
how to remove the factor $2$ in this inequality
at least for open manifolds $M$.

(ii) The energy-capacity inequality
$c_G (A) \le 2\, e (A,M)$ was proved in \cite{LM1} 
for every subset $A$ of {\it any}\, symplectic manifold 
$(M, \go )$.
This inequality implies that the Hofer norm on the group of
compactly supported Hamiltonian diffeomorphisms is non-degenerate.
In view of Corollary~1, it is conceivable that the factor $2$ can
always be omitted.
\diam
}
\end{remark}
Let $Z^{2n}(r)$ be the standard symplectic cylinder $B^2(r)
\times \RR^{2n-2} \subset \left( \RR^{2n}, \go_0 \right)$.
Another consequence of Theorem~1\,(ii) is
\begin{corollary}  \label{c:cap}
Assume that $\gs$ is an action selector for $\left( \RR^{2n},
\go_0 \right)$. Then for all $r > 0$,
\[
c_\gs \left( B^{2n}(r) \right) = c_\gs \left( Z^{2n}(r) \right)
= \pi r^2
\;\text{ and }\;
c^\gs \left( B^{2n}(r) \right) = c^\gs \left( Z^{2n}(r) \right)
= \pi r^2 .
\]
\end{corollary}

\subsection{An estimate for the smallest spectral value}

\ni
A {\it hypersurface} $S$ in $(M,\go)$
is a smooth compact connected orientable codimension
$1$ submanifold without boundary contained in $M \setminus \pp M$.
A {\it closed characteristic}\, on $S$ is an embedded circle in $S$ all
of whose tangent lines belong to the distinguished line bundle 
\[
\cl_S \,=\, \left\{ (x, \xi) \in TS \mid \go(\xi, \eta) =0 
\text{ for all } \eta \in T_x S \right\} .
\]
Denote by $\cp (S)$ the set of closed characteristics on $S$, 
and by $\cp^\circ (S)$ the set of those closed characteristics on $S$
which are contractible in $M$.
For $x \in \cp^\circ (S)$ let $\cd (x)$ be the set of smooth
discs $\bar x \colon D^2 \to M$ with boundary $x$.
The {\it contractible action spectrum}\, of $S$ is the set 
\[
\Sigma^\circ (S) \,=\, 
\left\{ \int_{\bar x} \go \;\Bigg|\; x \in \cp^\circ (S), \bar x \in
\cd (x)  \right\} .
\]
If $\go |_S = d \gl$ for some $1$-form $\gl$ on $S$, 
the {\it full action spectrum}\, is defined as 
\[
\Sigma (S,\gl) \,=\, 
\left\{ \int_x \gl \;\Bigg|\; x \in \cp (S) \right\} .
\]
Note that $\Sigma (S,\gl)$ is independent of the choice of $\gl$ if
$H^1(S;\RR) =0$.
The sets $\Sigma^\circ (S)$ and $\Sigma (S,\gl)$ are important
collections of numerical invariants of $S$,
see \cite{C,CFHW}.
Here, we are interested in the ``smallest'' spectral value.
If $\cp^\circ (S)$ or $\cp (S)$ is non-empty, we define
\[
\ga_1^\circ (S) = 
\inf \left\{ \left| \ga \right| \mid \ga \in \Sigma^\circ (S) \right\}
\quad \text{and} \quad
\ga_1 (S,\gl) \,=\, 
\inf \left\{ \left| \ga \right| \mid \ga \in \Sigma (S,\gl) \right\} .
\]
If $(M, \go)$ is not rational, then $\ga_1^\circ (S) =0$ for any
hypersurface $S \subset (M,\go)$ with $\cp^\circ (S) \neq
\emptyset$,
and our results for $\ga_1 \left( S, \gl \right)$ will deal with
hypersurfaces in weakly exact symplectic manifolds.
We shall thus assume in this paragraph that $(M,\go)$ is
rational.
We shall also assume that $S$
is a hypersurface of {\it contact type}. This means that
there exists a Liouville vector field $X$ 
which is defined near $S$ and is transverse to $S$.
Equivalently, there exists a contact form $\gl$ on $S$
(i.e., a $1$-form $\gl$ such that $d \gl = \go |_S$
and $\gl \wedge (d\gl)^{n-1}$ is a volume form on $S$). 
The equivalence is given by $\iota_X \go = \gl$.
A hypersurface $S$ is said to be of {\it restricted contact type}\, if
it is transverse to a Liouville vector field $X$ defined {\it on all}\,
of $M$.
The contact form $\gl = \iota_X \go$ is then globally defined
and $d \gl = d\iota_X \go = \go$ so that
$(M,\go)$ is exact.
If $(M,\go)$ is exact, then every hypersurface $S$ of contact
type with $H^1(S;\RR) =0$ is of restricted contact type.
\begin{examples}\
{\rm
{\bf 1.}
Consider a hypersurface $S$ of $\left( \RR^{2n}, \go_0 \right)$.
It bounds a bounded domain $U$. Then
\begin{eqnarray*}  
 \begin{array}{ll} 
             & U \text{ is convex}  \\
 \Longrightarrow  & U \text{ is starshaped} \\
 \Longrightarrow  & S \text{ is of restricted contact type} \\
 \Longrightarrow  & S \text{ is of contact type.}
 \end{array}
\end{eqnarray*}
Examples show that none of these arrows can be inverted, 
see \cite{B,He2}.

\s
\ni
{\bf 2.}
Let $T^*B$ be the cotangent bundle over a closed base endowed with the
symplectic form $\go_0 = d \gl_0$, where $\gl_0 = \sum_i p_i
dq_i$.

\s
(i)
Assume that $B$ is endowed with a Riemannian metric. Any regular
energy level $S_c = \{ H = c \}$ of a classical Hamiltonian
$H(q,p) = \frac 12 |p|^2 + V(q)$ is of contact type, 
see \cite[Theorem~1.2.2]{AH},
and if $c > \max V$, then $S_c$ is of restricted contact type,
see (iii) below.
If $c > \max V$, then $\Sigma (S_c) \neq \emptyset$, and if 
$c< \max V$, then $\Sigma^\circ (S_c) \neq \emptyset$; see
\cite[p.\ 131]{HZ} for a brief history of this existence problem
and for further references.

\s
(ii)
Assume that $S$ is a hypersurface of contact type in $\left(
T^*B, \go_0 \right)$ and that $\dim B \ge 2$.
Then $\Sigma (S) \neq \emptyset$ if the bounded component of
$T^*B \setminus S$ contains $B$, see \cite{HV}, and 
$\Sigma^\circ (S) \neq \emptyset$ if $B$ is simply connected,
\cite{V1}.

\s
(iii)
Assume that $S$ is a hypersurface in $T^*B$ such that for each $q \in B$
the intersection $S \cap T^*_q B$ bounds a strictly convex domain in
$T^*_qB$ containing $0$.
Then $S$ is transverse to the Liouville vector field $X(q,p) =
\sum_i p_i \frac{\pp}{\pp p_i}$, and hence $S$ is of restricted
contact type. Moreover,
there exists a unique Finsler metric $F \colon TB \to \RR$ on $ B$
such that $S$ is the unit cosphere bundle
$\left\{ x \in T^*B \mid F^* (x) =1 \right\}$,
where $F^* \colon T^*B \to \RR$ is the ``dual norm'' in each fibre,
defined by
\[
F^*(q,p) \,=\, \sup_{0 \neq \dot q \in T_qB} \frac{p(\dot q)}{F(q, \dot
q)} , 
\]
see \cite[Section 4.a)]{C}.
The projection $\pi \colon T^*B \to B$ induces a bijection between $\cp
(S)$ and the non-empty 
set of prime geodesics in the Finsler
metric $F$, and the action $\left| \int_x \gl_0 \right|$ of 
$x \in \cp (S)$ equals the $F$-length of $\pi (x)$. 
In particular, $\ga_1 (S, \gl_0)$ is the length of the shortest
closed $F$-geodesic on $B$.

\s
\ni
{\bf 3.}
Let $S$ be a smoothly embedded loop in the $2$-sphere $S^2$ endowed with
an area form. Then $S$ is of contact type, but not of restricted
contact type. It bounds two discs of areas $a_1$ and $a_2$, and
$\ga_1^\circ (S) = \min \left( a_1, a_2 \right)$.
}
\diam
\end{examples}

The following two propositions seem to be well known; 
proofs are given for the reader's convenience in
Appendix~\ref{app:2}.
\begin{proposition}  \label{p:spec}
Assume that $S$ is a hypersurface of contact type in a weakly
exact symplectic manifold $(M,\go)$.

\s
\ni
(i) 
If $S$ is simply connected or if $S$ is of restricted
contact type, then 
$0 \notin \Sigma^\circ (S)$ and $\Sigma^\circ (S)$ is closed.
If $\cp^\circ (S) \neq \emptyset$, we thus have $\ga_1^\circ (S) >0$. 

\s
\ni
(ii)
If $S$ is of contact type, then 
$0 \notin \Sigma (S, \gl )$ and $\Sigma (S, \gl)$
is closed for any contact form $\gl$ on $S$.
Thus, $\ga_1 (S,\gl) >0$ if $\cp (S) \neq \emptyset$. 
\end{proposition}

\begin{proposition}  \label{p:restricted}
Assume that $(M,\go)$ is a symplectic manifold, and that
$U \subset M$ is a relatively compact domain whose boundary $S$ is
a hypersurface of restricted contact type. 
If $\cp^\circ (S) \neq \emptyset$, then 
\[
0 \,<\, \ga^\circ_1 (S) \,\le\, c_{\HZ}^\circ (U,M) ,
\]
and if $\cp (S) \neq \emptyset$, then 
\[
0 \,<\, \ga_1 (S, \gl) \,\le\, c_{\HZ} (U) 
\]
for any globally defined contact form $\gl$.
\end{proposition}
For a convex bounded domain with smooth boundary in $\left(
\RR^{2n}, \go_0 \right)$, it holds that 
$\ga_1 (S) = c_{\HZ} (U)$, see \cite{HZ2}, and
for the full Bordeaux bottle (which is star--shaped) 
it holds that $\ga_1 (S) < c_{\HZ} (U)$, see \cite[p.\ 99]{HZ}.
\begin{problem}  \label{open:<}
{\rm
Assume that a hypersurface $S \subset B^4 (1) \subset \left(
\RR^{2n}, \go_0 \right)$ bounds a star-shaped domain $U$.
Then $c_{\HZ}(U) \le c_\gs(U)$ for any action selector $\gs$ for
$B^4(1)$ in view of Theorem~1\,(ii).
Does equality hold for $\gs_{V}$, $\gs_{\HZ}$ and $\gs_{\PSS}$?
Since it is known that $c_{\gs}(U)$ belongs to $\Sigma^\circ
(S)$ for every monotone action selector $\gs$, the first step toward
a solution of this problem would be to see whether $c_{\HZ} (U) \in
\Sigma^\circ (S)$.
\diam
}
\end{problem}

\m
Our next goal is to find upper bounds of $\ga_1^\circ (S)$ for
contact type hypersurfaces $S$ which might not be of restricted
contact type and which might not bound a domain $U$.
Since the interior of a hypersurface $S$ is empty, $c_\gs (S,M) =0$.
We thus consider the {\it outer spectral capacity} 
\[
\hat{c}_\gs \left( S,M \right) \,=\,
\inf \left\{ c_\gs \left( U,M \right) \mid S \subset U,\, U
\text{ open in } M \right\} .
\]
The {\it index of rationality} $\rho (M, \go) \in (0,\infty]$ of a rational
symplectic manifold $(M,\go)$ is set $\infty$ if $(M, \go)$ is
weakly exact and is defined to be the positive generator of the
cyclic subgroup $\go (\pi_2(M))$ of $\RR$ if $(M,\go)$ is not
weakly exact.
The following result improves Corollary~11.2 in \cite{FS}.
\\
\\
{\bf Theorem~2.}
{\it Assume that $(M,\go)$ is a rational symplectic manifold
admitting a weak action selector $\gs$, 
and consider a hypersurface $S \subset (M,\go)$ of contact type.
Then $\hat{c}_\gs (S,M) \le e(S,M)$, and if 
$\hat{c}_\gs \left( S,M \right) < \rho (M, \go)$, then
$\cp^\circ (S) \neq \emptyset$ and
\[
\ga_1^\circ (S) \,\le\, \hat{c}_\gs \left( S,M \right) \,\le\, e
\left( S,M \right) . 
\]
}

\ni
{\bf Corollary~2.}
{\it
Assume that an exhaustion of the rational symplectic manifold
$(M,\go)$ admits a weak action selector, 
and consider a hypersurface $S$ of contact type such that
$e \left( S,M \right) < \rho (M, \go)$. Then $\cp^\circ (S) \neq
\emptyset$ 
and
$\ga_1^\circ (S) \le  e \left( S,M \right)$.
}

\b
\begin{remarks}\
{\rm
{\bf 1.}
Again, Theorem~2 and Corollary~2 are sharp: For the circle
$S \subset \left( \RR^2, \go_0 \right)$ bounding $B^2(r)$ we
have $\ga_1^\circ(S) = e (S,\RR^2) = \pi r^2$. 

\s
\ni
{\bf 2.}
Theorem~2 and Corollary~2 establish 
the Weinstein conjecture 
for a large class of hypersurfaces of contact type. We refer to
\cite{Gi3, Schl} for the state of the art of the Weinstein conjecture.

\s
\ni
{\bf 3.}
Consider a contact hypersurface $S$ in $\left( \RR^{2n}, \go_0
\right)$.
Then $S$ is contained in a ball of radius $\Diam (S)$, where
$\Diam (S)$ is the diameter of $S$.
Since $e \left( B^{2n} (r),\RR^{2n} \right) = \pi r^2$, we find 
$e (S,\RR^{2n}) \le \pi \Diam (S)^2$, and hence
\[
\ga_1 (S) \,\le\, \pi \Diam (S)^2 ,
\] 
improving the estimates in \cite{FS,HZ1}.
}
\diam
\end{remarks}

We finally consider a billiard table $\overline{U} \subset
\RR^n$ with smooth boundary, i.e., $U$ is a bounded domain in
$\RR^n (q)$ with smooth connected boundary.
The length of a billiard trajectory on $\overline{U}$ is
measured with respect to the Euclidean length, and the Euclidean
volume is $\Vol (U)$.
Let $D^n$ be the closed unit ball in $\RR^n(p)$.
Combining a generalization of Corollary~2 with the construction
in \cite{V3}, we shall obtain the following result of Viterbo,
whose proof in \cite{V3} uses the fact that $\overline{U} \times D^n
\subset \left( \RR^{2n}, \go_0 \right)$
can be approximated by domains with boundary of restricted
contact type.
\\
\\
{\bf Proposition~2} (Viterbo){\bf .}
{\it
Let $U$ be a bounded domain in $\RR^n$ with smooth boundary. 
Then there exists a periodic billiard trajectory on
$\overline{U}$ of length
\[
l \,\le\, e \left( \overline{U} \times D^n, \RR^{2n} \right)
\,\le\, C_n \left( \Vol (U) \right)^{1/n} ,
\]
where $C_n$ is a constant depending only on $n$.
}

\b
\ni
{\bf Acknowledgements.}
We are grateful to Peter Albers, David Hermann and Matthias Schwarz 
for most helpful discussions.
We thank Dietmar Salamon for sending us Chapter~12 of the
forthcoming book \cite{MS3}.
Part of this note was written when the third author was visiting 
Hokkaido University; he wishes to thank Kaoru Ono for his warm hospitality.
The first and third author thank Hokkaido University and 
Leipzig University for their support.

\section{Proofs}
\label{sec:proofs}
\subsection{Proof of Theorem~1}
(i) Assume that $(M, \go)$ is a symplectic manifold admitting a
weak action selector $\gs$.
The support $\supp H$ of $H \in \ch$ is defined as the closure
of the set
\[
\bigcup_{t \in [0,1]} \left\{ x \in M \mid H_t(x)
\neq 0 \right\} .
\]
The main ingredient in the proof of Theorem~1 is
\begin{proposition}  \label{p:dis}
Assume that $H,K \in \ch$ are such that $\gf_K$ displaces $\supp H$.
Then $\gs (H) \le \left\| K \right\|$.
\end{proposition}
\proof
We follow the argument in \cite{FS,Gi3,HZ,Sch}.
Given a smooth function $\gl \colon [0,1] \ra [0,1]$ with
$\gl(0)=0$ and $\gl(1)=1$, we define $H^\gl \in \ch$ as
\[
H^{\gl} (t,x) \,=\, \gl'(t) \,H \left( \gl(t),x \right) .
\]
\begin{lemma}  \label{l:rep}
$\sigma \left( H^\gl \right) = \sigma (H)$.
\end{lemma}
\proof
The Hamiltonian flow of $H^\gl$ is a reparametrization of the
flow of $H$:
\[
\gf^t_{H^\gl}(x) \,=\, \gf^{\gl(t)}_H(x),  \quad\, x \in M ,
\]
Since $\gl(0)=0$ and $\gl(1)=1$, this reparametrization gives rise
to a one-to-one correspondence between one-periodic orbits of the flows:
the orbit $x^\gl \in\cp^\circ \left( H^\gl \right)$ corresponds to
$x \in \cp^\circ (H)$ when
\[
x^{\gl} (t) \,=\, x \left( \gl(t) \right) .
\]
Moreover,
\[
\int_0^1 H^\gl \left( t, x^\gl(t) \right) dt \,=\, \int_0^1 \gl'
(t)\, H \left( \gl(t),x \left(\gl(t) \right)\right) dt
\,=\, \int_0^1 H \left( t, x (t) \right) dt ,
\]
and hence $\ca_{H^\gl} \left( \bar{x}^\gl \right) = \ca_H \left(
\bar{x} \right)$ for all $(x, \bar x) \in \bar{\cp}^\circ (H)
\cong \bar{\cp}^\circ \left( H^\gl \right)$.
Consider now the smooth family of functions 
$\gl_\tau \colon [0,1] \ra [0,1]$ given by
\[
\gl_\tau (t) \,=\, (1-\tau) t + \tau \gl (t), \quad\, \tau \in
[0,1] .
\]
Then $H^{\gl_0} =H$ and $H^{\gl_1} = H^\gl$.
By the above, $\Sigma^\circ \left( H^{\gl_\tau} \right) =
\Sigma^\circ (H)$ for all $\tau \in [0,1]$, and according 
to \cite{HZ, Sch} and \cite[Lemma 2.2]{Oh3}, 
the Lebesgue measure of this set vanishes.
In view of (AS4), the map
\[
[0,1] \ra \Sigma^\circ(H), \quad\, \tau \mapsto 
\gs \big(  H^{\gl_\tau} \big)
\]
is constant. In particular, $\gs \left( H^\gl \right) = \gs (H)$.
\proofend

Utilizing Lemma~\ref{l:rep}, we can assume that $H_t=0$ for $t
\in [0,1/2]$ and $K_t=0$ for $t \in [1/2,1]$.
With this parametrization, 
\[
\left( \tau H \# K \right) (t,x) \,=\, \tau H (t,x) + K (t,x)
\] 
for all $\tau \in [0,1]$ and $(t,x) \in I \times M$, and the
time-$1$-flow of $\tau H \# K$ is the time-$1$-flow of $K$
followed by the time-$1$-flow of $\tau H$.
Since $\gf_K$ displaces $\supp \tau H = \supp H$ for each $\tau
\in (0,1]$, it follows that 
$x(0) \notin \supp \tau H$ for each
$x \in \cp^\circ (K)$ and each 
$x \in \cp^\circ \left( \tau H \# K \right)$, and we conclude
that  
$\cp^\circ \left( \tau H \# K \right) = \cp^\circ (K)$.
Therefore,
$\Sigma^\circ \left( \tau H \# K \right) = \Sigma^\circ (K)$ for each $\tau
\in [0,1]$, and arguing as above we find
$\gs \left( H \#K \right) = \gs (K)$.
The inverse $\gf_K^{-1}$ of $\gf_K$ is generated by
$K^- (t,x) = - K \left(t, \gf_K^t(x) \right)$.
Notice that $K \# K^- =0$. Combining this with (AS5), we conclude
\[
\gs(H) \,=\, \gs \left( H \# K \# K^- \right) \,\le\, \gs \left(
H \# K \right) + E^+ \left( K^- \right) \,=\, \gs(K) + E^+
\left( K^- \right) .
\]
Using
\[
\max_{x \in M} K^-(t,x) \,=\, \max_{x\in M} \left( -K \left(
t,\gf_K^t(x) \right) \right) \,=\, - \min_{x\in M} K (t,x)
\]
and (AS3) we finally obtain
\begin{eqnarray*}
\gs (H) \,\le\, \gs(K) + E^+ \left( K^- \right)
        &\le&  \int_0^1 \left( \max_{x\in M} K(t,x) -
                       \min_{x\in M} K(t,x) \right) \,dt \\
        &=& \left\| K \right\| , 
\end{eqnarray*}
as desired.
\proofend

Consider again a symplectic manifold $(M,\go)$ admitting a weak action
selector $\gs$, and let $A \subset M$.
In order to show that 
$c^\gs (A,M) \le e (A,M)$, 
we can assume that $e \left( A,M \right) <\infty$.
Fix $\gd >0$ and choose $K \in \ch$ such that $\left\| K
\right\| \le e \left( A,M \right) + \gd$ and $\gf_K$ displaces
$A$.
In particular, $\gf_K$ displaces $\supp H$ for any $H
\in \ch (I \times A)$.
Proposition~\ref{p:dis} thus yields $\gs (H) \le \left\| K
\right\| \le e \left( A,M \right) + \gd$. Taking the supremum
over $H \in \ch (I \times A)$, we find $c^\gs \left( A,M \right) \le e
\left( A,M \right) + \gd$, and since $\gd >0$ is arbitrary, the
claim follows.

\s
\ni
(ii)
Assume that $\gs$ is an action selector for $(M,\go)$. In order to show that
$c_{\HZ}^\circ \left( A,M \right) \le c_\gs (A,M)$,
we need to prove that $\max H \le \gs (H)$ for every $H \in
\cs_{\HZ}^\circ (A,M)$, and this holds by (AS2$^+$).
\proofend

\subsection{Proof of Theorem 2}
We consider an arbitrary hypersurface $S$ of a symplectic
manifold $(M,\go)$.
Examples show that $\Sigma^\circ (S)$ can be empty, see \cite{Gi, GG}.
We therefore follow \cite{HZ1} and consider parametrized neighborhoods of $S$.
Since $S$ is orientable and contained in $M \setminus \pp M$, 
there exists an open neighborhood $I$ of $0$ and a smooth diffeomorphism 
\[
\psi \colon S \times I \,\ra\, U \subset M
\]
such that $\psi (x,0) =x$ for $x \in S$.
We call $\psi$ a {\it thickening of $S$} and 
set $S_\eps = \psi \left( S \times \left\{ \eps \right\}
\right)$.

\begin{theorem}  \label{t:action}
Assume that $(M,\go)$ is a rational symplectic manifold
admitting a weak action selector $\gs$, and consider a
hypersurface $S$ of $(M,\go)$.
Then $\hat{c}_\gs (S,M) \le e (S,M)$.
If $\hat{c}_\gs (S,M) < \rho (M,\go)$, 
then for every thickening $\psi$ of $S$ and 
every $\gd > 0$ there exists $\eps \in \left[ -\gd, \gd \right]$
such that 
\[
\cp^\circ \left( S_\eps \right) \neq \emptyset 
\quad \text{and} \quad
\ga_1 \left( S_\eps \right) \le \hat{c}_\gs \left( S,M \right) + \gd 
                            \le e \left( S,M \right) + \gd .
\]
\end{theorem}

\proof
By Theorem~1\,(i), $c_\gs (U,M) \le e (U,M)$ for every
open neighborhood $U$ of $S$, and hence the first claim follows
from taking the infimum. 
In order to prove the second claim,
we can assume that $\gd > 0$ is so small that 
\[
C \,:=\, \hat{c}_\gs (S,M) +\gd \,<\, \rho := \rho (M,\go) .
\] 
Let $\tau \in (0,\gd)$ be so small that for the neighborhood 
$U_\tau = \psi \left( S \times (-\tau, \tau) \right)$ of $S$ we
have
$c_\gs \left( U_\tau, M \right) < C$.
We choose a smooth function $f \colon \RR \ra [0,C]$ such that
\begin{eqnarray*}  
 \begin{array}{rcll} 
    f(t) &=& 0 & \text{ if $t \notin 
            \, \left]-\tfrac{\tau}{2}, \tfrac{\tau}{2} \right[$}, \\[0.4em]
    f(t) &=& C & \text{ if $t \in \left[ -
\tfrac{\tau}{4},\tfrac{\tau}{4} \right]$},  
                                                                  \\[0.4em]
    | f'(t) | &\neq& 0  & \text{ if $t \in \,
               \left]-\tfrac{\tau}{2}, -\tfrac{\tau}{4} \right[
          \cup \left] \tfrac{\tau}{4}, \tfrac{\tau}{2} \right[$}.  
 \end{array}
\end{eqnarray*}
Define $H \in \cs (M)$ by
\begin{equation*}  
H (x) \,=\, 
 \left\{
  \begin{array}{ll}
   f(t) & \text{if } x \in S_t , \\ [0.2em]
   0    & \text{otherwise}.
  \end{array}
 \right.
\end{equation*}
By (AS2), we have $\gs (H) >0$, so that
\begin{equation}  \label{e:wein:0ge}
0 \,<\, \gs (H) \,\le\, c_\gs \left( U_\tau,M \right) \,<\, C
\,<\, \rho .
\end{equation}
By (AS1), there exists 
$(x,\bar{x}) \in \bar\cp^\circ (H)$ such that 
$\gs (H) = \ca_H \left( \bar x \right)$.
\begin{lemma}  \label{l:nonconstant}
The orbit $x$ is not constant.
\end{lemma}
\proof
If $x$ is constant, our choice of $H$ yields $H(x) \in \left\{
0,C \right\}$. If $H(x)=0$, then
\begin{equation}  \label{e:0}
0 \,<\, \gs (H) \,=\, \ca_H \left( \bar x \right) \,=\, -\int_{\bar
x } \go  \,<\, C ,
\end{equation}      
and if $H(x)=C$, then

\begin{equation}  \label{e:C}
0 \,<\, \gs (H) \,=\, \ca_H \left( \bar x \right) \,=\, -\int_{\bar
x } \go + C \,<\, C .
\end{equation}
Since $\bar x$ is a sphere, $\int_{\bar x }\go \in \rho \ZZ$, and
so both \eqref{e:0} and \eqref{e:C} contradict $C< \rho$.
\proofend

By construction of $H$ and by Lemma~\ref{l:nonconstant}, 
there exists $\eps \in (-\tau, \tau) \subset [-\gd, \gd]$ such that 
$(x,\bar x) \in \bar \cp^\circ \left( S_\eps \right)$ and $0 <
f(\eps) < C$ and
\[
0 \,<\, \gs (H) = -\int_{\bar x} \go + f(\eps) \,<\, C .
\]
We conclude that $\left| \int_{\bar x} \go \right| \le C$, so
that $\ga_1 \left( S_\eps \right) < \hat{c}_\gs \left( S,M \right) +
\gd$, as claimed.
\proofend

Theorem~2 is an easy consequence of Theorem~\ref{t:action};
we refer the reader to \cite[Section~11]{FS} for a detailed argument.

\begin{remark}
{\rm
Assume that $S$ is a hypersurface of contact type bounding a domain
$U$. An obvious modification of the proof of Theorem~\ref{t:action}
shows that Theorem~2 holds for $c_\gs \left(U,M \right)$ as
well.
It is, however, easy to see that $\hat{c}_\gs \left( S,M \right)
\le c_\gs \left( U,M \right)$.
}
\end{remark}

\subsection{Proof of Proposition~2}
Let $U \subset \RR^n$ be a bounded domain with smooth boundary,
and abbreviate $e = e ( \overline{U} \times D^n, \RR^{2n})$.
We choose a smooth family $W_t$, $t \in (0,1]$, of bounded
domains in $\RR^{2n}$ with smooth boundaries $S_t$ such that
\[
W_t \subset W_{t'} \;\text{ if }\, t \le t'
\quad \text{ and } \quad
\bigcap_{t \in (0,1]} W_t = \overline{U} \times D^n .
\]
In view of Theorem~\ref{t:action}, we can find a sequence $t_k \ra 0$
such that
$S_{t_k}$ carries a closed characteristic $\gg_k$ with
$\int_{\gg_k} \gl_0 \le e + 1/k$.
Choose a parametrization $\gg_k (t) = \left( q_k(t), p_k(t)
\right) \colon [0,1] \ra S_{t_k}$ of $\gg_k$.
Since $p_k(t) = f_k(t) \frac{\dot{q}_k (t)}{\left| \dot{q}_k(t)
\right|}$ for a function $f_k \ge 1$,
\[
\length (q_k) \,=\, \int_0^1 \left| \dot{q}_k (t) \right| dt
\,\le\, \int_0^1 \langle p_k(t), \dot{q}_k(t) \rangle dt
\,=\, \int_{\gg_k} \gl_0 \,\le\, e + \tfrac 1k.
\]
We therefore find a subsequence $q_{k_m}$ converging to a
billiard trajectory $q$ on $U$, see \cite{BG}, and $\length (q)
\le e$. 
The estimate $e \le C_n \left( \Vol (U) \right)^{1/n}$ is proved
in \cite{V3}.
\proofend

\appendix

\section{Constructions of action selectors}  \label{app:0}
\ni
In this appendix we outline the constructions of the (weak)
action selectors $\gs_{\V}$, $\gs_{\HZ}$ and $\gs_{\PSS}$.

\s
\ni
{\bf 1. $\gs_{\V}$ and $\gs_{\HZ}$ for $\left( \RR^{2n}, \go_0 \right)$.}

\s
\ni
{\bf The selector $\gs_{\V}$.}
The following description is taken from \cite{H3,V}.
The construction of the action selector $\gs_{\V}$ makes use of
generating functions.
For $H \in \ch \left( \RR^{2n} \right)$ the graph
\[
\Gamma_H \,=\, \left\{ (x, \gf_H(x) \mid x \in \RR^{2n} \right\} 
\]
of $\gf_H$ is a Lagrangian submanifold of $\left( \RR^{2n} \times
\RR^{2n}, -\go_0 \oplus \go_0 \right)$.
Under the symplectomorphism
\[
\left( q,p,Q,P \right) \,\mapsto\, 
\left( \frac{q+Q}{2}, \frac{p+P}{2}, p-P,Q-q \right)
\]
of $\left( \RR^{2n} \times \RR^{2n}, -\go_0 \oplus \go_0 \right)$ with the cotangent bundle $T^*
\Delta$ over the diagonal in $\RR^{2n} \times \RR^{2n}$ the graph
$\Gamma_H$ is mapped to a Lagrangian submanifold of $T^*\Delta$ which
outside a compact set coincides with $\Delta$. Adding a point $\infty$
to $\Delta \cong \RR^{2n} \cong S^{2n} \setminus \{ \infty \}$ one
obtains a closed Lagrangian submanifold $\Gamma$ of $T^*S^{2n}$, which is
Hamiltonian isotopic to the zero section.
According to a result of Sikorav, $\Gamma$ admits a generating function
quadratic at infinity. This means that for some $N \ge 0$ there exists a
smooth function $S \colon S^{2n} \times \RR^N \ra \RR$ such that
\begin{itemize}
\item[$\bullet$]
$0$ is a regular value of the fibre derivative $\pp_\xi S \colon S^{2n}
\times \RR^N \ra \RR^N$;
\item[$\bullet$]
$\Gamma = 
\left\{ \left( q, \pp_q S(q, \xi) \right) \mid (q,\xi) \in S^{2n}
\times \RR^N \text{ satisfies } \pp_\xi S(q,\xi) =0 \right\}$;
\item[$\bullet$]
$S(q,\xi) = Q(\xi)$ away from a compact set for a non-degenerate
quadratic form $Q$ on $\RR^N$.
\end{itemize}
Critical points $(q,\xi)$ of $S$ correspond to $1$-periodic orbits $x \in
\cp^\circ (H)$, and if 
$S$ is normalized such that $S(\infty,0) =0$, then $\ca_H(x) =
-S(q,\xi)$. For $a \in \RR$ set
\[
S^a \,=\, \left\{ (q,\xi) \in S^{2n} \times \RR^N \mid S(q,\xi) < a \right\} .
\]
We denote by $\H_*$ homology with real coefficients.
It is shown in \cite{Th,V}, see also \cite{Tr},  
that the relative homology groups $\H_* \left( S^b, S^a \right)$ do not
depend on the choice of $S$, and it is easily seen that $\H_* \left( S^a,
S^{-a} \right)$ does not depend on $a$ for $a$ large enough.
Let $i$ be the index of $Q$, and for $a$ large enough let $1 \in \H_i
\left( S^a, S^{-a} \right)$ be the image of the positive generator of
$\H_0 \left( S^{2n} \right)$ under the Thom isomorphism
$\H_0 \left( S^{2n} \right) = \H_i \left( S^a, S^{-a} \right)$.
Consider the map
\[
j_S^\gl \colon \H_i \left( S^a, S^{-a} \right) \,\ra\, \H_i \left( S^a, S^\gl \right)
\]
induced by inclusion. The action selector $\gs_{\V}$ is defined as
\[
\gs_{\V} (H) \,=\, - \inf \left\{ \gl \mid j_S^\gl (1) =0 \right\} .
\]
It follows from classical Lusternik--Schnirelmann theory that $\gs_{\V}
(H)$ is a critical value of $-S$ and hence belongs to $\Sigma^\circ (H)$.
Similarly, the axioms (AS2) to (AS5) are translations of
properties of special critical values of generating functions
quadratic at infinity, which can be verified by elementary but by no
means trivial topological arguments.

\m
\ni
{\bf The selector $\gs_{\HZ}$.}
The following description is taken from Section~5.3 of \cite{HZ}.
The action selector $\gs_{\HZ}$ is defined via a direct minimax on the
space of loops 
\[
E \,=\, H^{1/2} \,=\, \left\{ x \in L^2 \left( S^1; \RR^{2n} \right)
\,\Bigg|\, \sum_{k \in \ZZ}  \left| k \right| \left| x_k \right|^2 < \infty
\right\} 
\]
where 
\[
x \,=\, \sum_{k \in \ZZ} e^{k 2\pi J t} x_k, \qquad x_k \in \RR^{2n} ,
\]
is the Fourier series of $x$ and $-J$ is the standard almost complex
structure of $\RR^{2n} \cong \CC^n$.
The space $E$ is a Hilbert space with inner product
\[
\langle x,y \rangle \,=\, \langle x_0,y_0 \rangle + 2\pi \sum_{k \in
\ZZ} \left| k \right| \langle x_k,y_k \rangle ,
\]
and there is an orthogonal splitting 
\[
E \,=\, E^- \oplus E^0 \oplus E^+, \qquad x = x^- + x^0 + x^+ ,
\]
into the spaces of $x \in E$ having only Fourier coefficients for $j<0$,
$j=0$, $j>0$.
The action functional $\ca_H \colon \cp^\circ (H) \ra \RR$
extends to $E$ as
\[
\ca_H(x) \,=\, a(x) + b(x), \qquad x \in E,
\]
where $a(x) = - \int_{\bar{x}} \go = 
\frac 12 \left\| x^+ \right\|^2 - \frac 12 \left\| x^-
\right\|^2$ 
and $b(x) = \int_0^1 H(t,x(t)) \,dt$.
The function $\ca_H \colon E \ra \RR$ is differentiable, and its
gradient is given by
\[
\nabla \ca_H (x) \,=\, x^+ - x^- + \nabla b(x) .
\]
Prompted by the structure of the gradient flows generated by such gradients,
we consider the group $G$ of homeomorphisms $h$ of $E$ 
satisfying
\begin{itemize}
\item[$\bullet$]
$h$ and $h^{-1}$ map bounded sets to bounded sets;
\item[$\bullet$]
there exist continuous maps $\gg^\pm \colon E \ra \RR$ and $k
\colon E \ra E$ mapping bounded sets to precompact sets and such
that there exists $r=r(h)$ satisfying 
\[
\gg^\pm (x) =0
\quad \text{and} \quad k(x)=0 \quad\, \text{ for all }\; x \in E
\text{ with } \left\| x \right\| \ge r
\]
and
\[
h(x) \,=\, e^{\gg^+(x)} x^+ + x^0 + e^{\gg^- (x)} x^- + k(x) 
   \quad\, \text{ for all }\; x \in E .
\]
\end{itemize}
The selector $\gs_{\HZ}$ is defined as
\[
\gs_{\HZ} (H) \,=\, \sup_{F\in\cf} \inf_{x \in F} \ca_H(x)
\]
where 
\[
\cf \,=\, \left\{ h \left( E^+ \right) \mid h \in G \right\} .
\]
It follows from Schauder's fixed point theorem that
\[
h \left( E^+ \right) \cap \left( E^- \oplus E^0 \right) \neq
\emptyset \quad\, \text{ for all }\; h \in G .
\]
This intersection result implies that $\gs_{\HZ} (H)$ is finite.
Using the fact that $\ca_H$ satisfies the Palais--Smale
condition one then shows that $\gs_{\HZ} (H) \in \Sigma^\circ (H)$.
The other axioms readily follow from definitions.

\m
\ni
{\bf 2. $\gs_{\PSS}$ for weakly exact closed symplectic manifolds.}
The following description closely follows \cite{Sch}.
The construction of the action selector $\gs_{\PSS}$ is along
the lines of the construction of $\gs_{\V}$, but instead of working
with a finite dimensional reduction $S$, one works with the
action functional $\ca_H$ on the full space $\cl^\circ$ of
contractible loops.
Notice that $\cp^\circ (H)$ is the set of critical points of $\ca_H$.
Let $H \in \ch$ be such that all $x \in \cp^\circ (H)$ are
non-degenerate in the sense that
\begin{equation}  \label{e:det}
\det \left( \id - d \gf^1_{H}(x(0) \right) \neq 0.
\end{equation}
Then $\cp^\circ (H)$ is a finite set, and $\ca_H$ is a Morse
function on $\cl^\circ$. The Morse indices are, however, all
infinite and hence do not lead to a change of topology of the
sublevel sets of $\ca_H$. Floer overcame this problem by
constructing a {\it relative}\, Morse theory for $\ca_H$, 
from which he extracted a Morse-type homology isomorphic to the homology of
$M$. It is called Floer homology and will serve as a substitute for the
homology groups $\H_* \left( S^a, S^{-a} \right)$ with $a$ large.

Consider the $\ZZ_2$-vector space $\CF (M,H)$ freely generated by
the elements of $\cp^\circ(H)$. In order to define a
differential on it, fix $x,y \in \cp^\circ(H)$. For a generic 
family $J_t$, $t \in S^1$, of almost complex structures on $M$
for which $g_t = \go \circ \left( \id \times J_t \right)$ is a Riemannian
metric for each $t$, 
the set $\cm (x,y)$ of solutions $u \in C^\infty \left( \RR \times S^1, M
\right)$ of the elliptic partial differential equation 
\begin{equation}  \label{e:floer} 
   \partial_s u+J_t(u)(\partial_t u-X_{H_t}(u)) \,=\, 0
\end{equation}
with asymptotic boundary conditions
\begin{equation}  \label{e:asym} 
    \displaystyle\lim_{s \to -\infty}u(s,t)  = x(t), \qquad
    \displaystyle\lim_{s \to \infty}u(s,t) = y(t) ,
\end{equation}
is a smooth manifold. 
Notice that $u \in \cm(x,y)$ is a downward gradient flow line
of $\ca_H$ with respect to the $L^2$-metric on $\cl^\circ$
induced by the family $g_t$.
Let $\cm^1(x,y)$ be the union of the
$1$-dimensional components of $\cm (x,y)$. 
The real numbers act on $\cm^1(x,y)$ by shift in the $s$-variable.
Since $(M,\go)$ is weakly exact, 
$\cm^1(x,y) / \RR$ is compact, and so one can set
\[
n(x,y) \,=\,  \# \left\{ \cm^1 (x,y) / \RR \right\} \mod 2.
\]
The linear map $\pp$ on $\CF(M,H)$ defined as the linear
extension of
\begin{equation}  \label{e:pp}
\pp x \,= \sum_{y \in \cp^\circ (H)} n(x,y) \,y, 
\qquad x \in \cp^\circ (H),
\end{equation}
indeed satisfies $\pp \circ \pp =0$, and the homology $\HF(M,H)$
of the resulting complex is
independent of the choice of the family $J_t$.
It is in fact independent of $H$ as well, and agrees with the
homology $\H(M;\ZZ_2)$ of $M$.

One way of constructing an isomorphism from the homology
$\H(M;\ZZ_2)$ to the Floer homology $\HF(M,H)$
is as follows.
The homology $\H(M;\ZZ_2)$ is canonically isomorphic to the Morse
homology $\HM(M;\ZZ_2)$ of $M$. A chain complex $\CM(M,F)$ for
this homology is generated by the critical points of a 
Morse function $F$ on $M$, and the differential is as in
\eqref{e:pp}, where now $n(c,c')$ is the number (mod $2$)
of negative gradient flow lines of $F$ (with respect to a
generic Riemannian metric) connecting critical points $c,c' \in
\Crit F$.
In order to relate $\HM(M;\ZZ_2)$ with $\HF(M,H)$, one defines a
chain map
$\phi \colon \CM(M,F) \ra \CF(M,H)$ by
\[
\phi(c) = \sum_{x \in \cp^\circ(H)} n(c,x) \,x, \quad \, c \in \Crit(F) ,
\]
where now $n(c,x)$ is the number (mod $2$) of mixed trajectories
$(\gg,u)$ from $c \in \Crit F$ to $x \in \cp^\circ (H)$ as in
Figure~\ref{figure2}.

\begin{figure}[h] 
 \begin{center}
   \psfrag{c}{$c$}
   \psfrag{g}{$\gg$}
   \psfrag{gu}{$\gg (0)$}
   \psfrag{u}{$u$}
   \psfrag{x}{$x$}
  \leavevmode\epsfbox{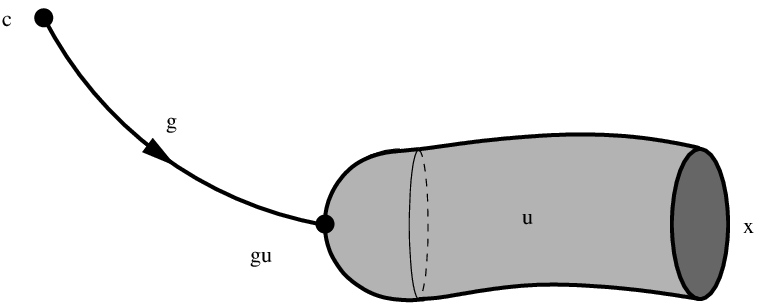}
 \end{center}
 \caption{An element of $\cm (c,x)$.}
 \label{figure2}
\end{figure}
%
%

\ni
Here, $\gg \colon (-\infty,0] \ra M$ is an integral curve of $ - \nabla F$,
i.e.,
\[
\dot{\gg}(s) \,=\, - \nabla F(\gg(s)), 
\]
and $u \colon \RR \times S^1 \ra M$ is a solution of the equation
\[
\pp_s u+J_t(u)(\pp_t u-X_{\gb(s) H_t}(u)) \,=\, 0 ,
\] 
where $\gb \colon \RR \ra [0,1]$ is a smooth cut off function such that 
\[
\gb (s) =0,\,\, s \le 0;
\quad\,
\gb'(s) \ge 0,\,\, s \in \RR;
\quad\,
\gb (s) =1,\,\, s \ge 1 ,
\]
and the boundary conditions are
\[
\lim_{s \ra -\infty} \gg(s) =c, 
\quad
\gg (0) = \lim_{s \ra -\infty} u(s,t),
\quad
\lim_{s \ra \infty} u(s,t) = x (t) .
\]
The map $\phi$ indeed commutes with the differentials and thus
induces a map $\Phi^H \colon \HM(M;\ZZ_2) \ra \HF (M,H)$, which
turns out to be an isomorphism. The composition
\[
\Phi_{\PSS}^H \colon \H(M;\ZZ_2) \cong \HM (M;\ZZ_2) 
\xrightarrow{\Phi^H} \HF(M,H)
\]
is called the Piunikhin--Salamon--Schwarz (PSS) isomorphism.

The Floer chain complex $\CF(M,H)$ comes with a natural real filtration
given by the action functional:
For $\gl \in \RR$ consider the linear subspace 
\[
\CF^\gl(M,H) \,=\, \left\{ \sum_{x \in \cp^\circ (H)} \xi_x \,x \,\Bigg|\,
\xi_x =0 \;\text{ if }\, \ca_H(x) >\gl \right\} 
\]
of $\CF(M,H)$.
Since $\ca_H$ decreases along solutions of \eqref{e:floer}, the
boundary operator $\pp$ preserves $\CF^\gl(M,H)$ and hence descends
to the quotient complex $\CF(M,H) / \CF^\gl(M,H)$.
Its homology is denoted by $\HF^\gl(M,H)$.
Let $j^\gl \colon \HF(M,H) \ra \HF^\gl(M,H)$ be the map induced by the
natural projection $\CF(M,H) \ra \CF(M,H)/ \CF^\gl(M,H)$, 
and let $1$ be the generator of $\H_{2n}(M;\ZZ_2) \cong \H^0(M;\ZZ_2)$.
Since $\Phi_{\PSS}^H$ is an isomorphism, 
\begin{equation}  \label{def:PSS}
\gs_{\PSS} (H) \,=\, 
  \inf \left\{ \gl \in \RR \mid j^\gl \left( \Phi_{\PSS}^H(1) \right) =
                                                0 \right\} 
\end{equation}
is finite. Furthermore, it is clear that $ \gs_{\PSS} (H)$ lies 
in $\Sigma^\circ (H)$.
So far we have defined $\gs_{\PSS}(H)$ for those $H \in \ch$
that satisfy \eqref{e:det}. 
The function $\gs_{\PSS}$ is $C^0$-continuous on this set;
this is proved by using the compatibility 
$\Phi_{\PSS}^K = \Phi^{KH} \circ \Phi_{\PSS}^H$ of the
$\PSS$-isomorphisms with the canonical isomorphism 
\[
\Phi^{KH} \colon \HF(M,H) \ra \HF(M,K)
\]
obtained by counting solutions of \eqref{e:floer} 
with $H$ replaced by $(1-\gb) H + \gb K$ and with asymptotic boundary
conditions $x \in \cp^\circ (H)$, $y \in \cp^\circ(K)$.
Since the set of $H$ satisfying \eqref{e:det} is 
$C^\infty$-dense in $\ch$, the function $\gs_{\PSS}$ 
$C^0$-continuously extends to all of $\ch$, and one readily
verifies that $\gs_{\PSS} (H) \in \Sigma^\circ (H)$ for all $H
\in \ch$. 
The explicit nature of the isomorphism $\Phi_{\PSS}$ is also important for
the verification of the remaining axioms.
The proof of (AS2) is difficult. One studies equation~\eqref{e:floer}
with $H$ replaced by $\gl H$ for each $\gl \in [0,1]$ and uses a
cobordism argument similar to the one used in proving $\pp \circ
\pp =0$. 
Finally, the proof of (AS5) uses the product structure on Floer
homology given by the pair of pants product and a sharp energy
estimate for the pair of pants.

\m
\ni
{\bf 3. $\gs_{\PSS}$ for weakly exact convex symplectic manifolds.}
The construction of $\gs_{\PSS}$ for weakly exact convex symplectic
manifolds in \cite{FS} follows the construction outlined above.
We can assume that $M$ is compact.
Since $M$ has boundary $\pp M$, the set $\cp^\circ(H)$ contains
infinitely many critical points for every $H \in \ch$.
This problem is overcome by using the geometry of $M$ near $\pp M$. 
On a neighborhood $\pp M \times (-\eps, 0]$ with coordinates
$(x,r)$ the symplectic form $\go$ is $d \left( e^r \ga \right)$, where
$\ga = ( \iota_X \go )|_{\pp M}$. 
Let $\ch_\pp$ be the set of Hamiltonians on $M$ which satisfy
\eqref{e:det} and which near $\pp M$ are of
the form $H (x,r) = f(r)$
with $f'(r)$ positive and so small that the flow
$\gf_H^t$ has no $1$-periodic orbits near $\pp M$.
For $H \in \ch_{\pp}$ the set $\cp^\circ (H)$ is
finite, and the strong maximum principle for uniformly elliptic
operators implies that solutions of \eqref{e:floer} and
\eqref{e:asym} stay away from $\pp M$. 
The Floer homology $\HF (M,H)$ can thus be defined and, by using a
yet stronger version of the maximum principle, the isomorphism $\Phi^H$
from the Morse homology $\HM (M;\ZZ_2)$ to $\HF (M,H)$ can also be
constructed. Here, the Morse homology is defined via Morse
functions which are of the form $F(x,r) = e^{-r}$ near $\pp M$.
After identifying $\HM_{2n} (M;\ZZ_2)$ with 
$\H_{2n} (M,\pp M;\ZZ_2) \cong \H^0 (M;\ZZ_2)$,
one defines $\gs_{\PSS} (H) \in \Sigma^\circ(H)$
as in \eqref{def:PSS}.
Since $\ch$ lies in the $C^\infty$-closure of $\ch_\pp$ and
since $\gs_{\PSS}$ is $C^0$-continuous on $\ch_\pp$, we
obtain a $C^0$-continuous function $\gs_{\PSS}$ on $\ch$
verifying axioms (AS1) and (AS2) for an action selector.
The remaining axioms are also verified as in the closed case.

Consider now an exhaustion $M_1 \subset M_2 \subset \cdots$ of a
non-compact convex weakly exact symplectic manifold $(M, \go)$,
and denote the action selector $\gs_{\PSS}$ for $M_i$ by
$\gs_i$.
The following examples shed some light on the problem whether an
action selector for an exhaustion of $(M,\go)$ fits together to
an action selector of $(M,\go)$.

\s
\ni 
(i) 
Assume that $\left( M \setminus M_1, \go \right)$ is
symplectomorphic to $\left( \pp M_1 \times (0,r_\infty), d
\left( e^r \ga \right) \right)$ for some $r_\infty \in (0,\infty]$
and that $M_i = M_1 \cup \pp M_1 \times [0,r_i]$, $0=r_1<r_2<
\dots < r_\infty$.
Then every $H_i \in \ch_{\pp} (M_i)$ extends to 
$H_{i+1} \in \ch_{\pp} (M_{i+1} )$ 
with $H_{i+1} (x,r) = f(r)$ on $M_{i+1}
\setminus M_i$.
The chain complexes $\CF (M_i, H_i)$ and $\CF (M_{i+1},
H_{i+1})$ then agree, and so  
$\gs_{i+1}$ restricts to $\gs_i$ on $M_i$.
Examples are $\left( \RR^{2n}, \go_0 \right)$ exhausted by
balls, cotangent bundles $\left( T^*B, \go_0 \right)$ over a
closed base exhausted by ball bundles, and, more generally,
Stein manifolds $(M,J,f)$. In the latter example, $(M,J)$ is an open 
complex manifold and $f \colon M \ra \RR$ is a smooth exhausting 
plurisubharmonic function  without critical points off a
compact subset of $M$.
Here, ``exhausting'' means that $f$ is proper and bounded from below,
and ``plurisubharmonic'' means that $\go = - d \left( df \circ J
\right)$ is a symplectic form with $\go \left( v,Jv \right) >0$
for all $0 \neq v \in TM$. 
Then the gradient $X= \nabla f$ with respect to the K\"ahler metric
$\go \circ \left( \id \times J \right)$ satisfies 
$\cl_X \go =\ \go$. Furthermore, if the critical points of $f$ are contained 
in, say, $\{ f < 1 \}$, the manifold $M$ is exhausted by the exact convex
symplectic manifolds $M_i = \{ f \le i \}$.

\s
\ni
(ii)
Assume that $M_{i+1} \setminus M_i$ is diffeomorphic but not
symplectomorphic to $\pp M_i \times \left( r_i,r_{i+1} \right]$.
Then every $H_i \in \ch_{\pp} (M_i)$ extends to
$H_{i+1} \in \ch_{\pp} (M_{i+1})$ such that the vector
spaces $\CF (M_i, H_i)$ and $\CF (M_{i+1}, H_{i+1})$ agree.
Moreover, then the homology groups of $M_i$ and $M_{i+1}$ agree, so that
the Floer homology groups $\HF \left(M_i, H_i \right)$ and 
$\HF \left(M_{i+1}, H_{i+1} \right)$ also agree.
However, solutions of \eqref{e:floer} and \eqref{e:asym} for
$H_{i+1}$ used to define the differential $\pp_{i+1}$ of
$\CF (M_{i+1}, H_{i+1})$ might enter $M_{i+1}
\setminus M_i$, and so the differentials $\pp_i$ and $\pp_{i+1}$
might not agree. It is thus unclear whether $\gs_i \left( H_i
\right) = \gs_{i+1} \left(H_{i+1}\right)$, and so $\gs_{i+1}$
might not restrict to $\gs_i$ on $M_i$.

\s
\ni
An example of a symplectic manifold admitting an exhaustion of type~(ii)
but no exhaustion of type (i) is the camel space $M \subset
\left( \RR^4, \go_0 \right)$ defined as
\[
M \,=\, \left\{ y_1 <0 \right\} \cup \left\{ y_1 >0 \right\}
\cup \left\{ x_1^2 + x_2^2 +y_2^2 < 1,\, y_1=0 \right\} ,
\]       
cf.\ \cite[Proposition~3.4.A]{EG}.
It is not hard to find an exhaustion of type~(ii): 
For $i \ge 1$ consider the union of the two closed $4$-balls of radius
$i$ centred at $y_1= \pm i$.
Smoothing this set appropriately near $y_1=0$ for each $i$ one obtains 
starshaped dumbbells $M_i$ with smooth boundary 
forming an exhaustion of
$M$ of type~(ii), see \cite[Lemma~5.1]{MT} for details.
Assume now that $M$ admits an exhaustion $(M_i)$ such that
$\left( M \setminus M_1, \go_0 \right)$ is symplectomorphic to
$\left( \pp M_1 \times (0,r_\infty), d (e^r \ga) \right)$.
Since $M$ has infinite volume, $r_\infty = \infty$.
Since $M = M_1 \cup \pp M_1 \times (0,\infty)$ is diffeomorphic to
the standard $\RR^4$, the interior of $M_1$ is also
diffeomorphic to the standard $\RR^4$, and so $M \setminus M_1$
is diffeomorphic to $S^3 \times \RR$. In particular, $H^1 \left(
\pp M_1; \RR \right) =0$.
The primitive $\gl = \iota_{\frac{\pp}{\pp r}} \go_0$ of $\go_0$
on $M \setminus M_1$ therefore smoothly extends to all of $M$.
The Liouville vector field $X$ defined by $\iota_X \go_0 = \gl$
then agrees with $\frac{\pp}{\pp r}$ on $M \setminus M_1$ and
hence integrates to a flow on $M$. A standard construction now
shows that the identical embedding of a dumbbell as above
extends to a symplectic embedding of $\left( \RR^4, \go_0
\right)$ into $M$, but this contradicts the Symplectic Camel
Theorem, see Lemma~4.1.7 and Proposition~4.1.8 in \cite{T}.
The same arguments apply to camel spaces in $\left( \RR^{2n}, \go_0
\right)$ for $n \ge 3$, where one then uses the Symplectic Camel
Theorem proved in \cite{V}.

\s
\ni
(iii)
If $M_{i+1} \setminus M_i$ is not diffeomorphic to $\pp M_i
\times (0,1]$, then any extension 
$H_{i+1} \in \ch_{\pp} (M_{i+1})$ of $H_i \in \ch_{\pp} (M_i)$ must
have critical points on $M_{i+1} \setminus M_i$. Already the vector
spaces $\CF (M_i, H_i)$ and $\CF (M_{i+1}, H_{i+1})$ are thus different, and so it is
unclear whether $\gs_i \left( H_i \right) = \gs_{i+1}
\left(H_{i+1}\right)$.
An example is the orientable surface $M$ of infinite genus
and with one end. More explicitly, $M$ is the surface with exhaustion $M_1
\subset M_2 \subset \cdots$, where $M_1$ is a torus with an open
disc removed and $M_{i+1} \setminus \Int M_i$ is a torus with two
open discs removed, $i = 1,2,\dots$.
For every area form on $M$ we obtain an exact
convex symplectic manifold, and it is unclear whether $\gs_{i+1}$
restricts to $\gs_i$ on $M_i$, $i = 1,2, \dots$. 
Higher-dimensional examples are products of such surfaces and also
Stein manifolds of infinite topology.

\m
\ni
{\bf 4. $\gs_{\PSS}$ for rational strongly semi-positive closed symplectic
manifolds.}
The construction of $\gs_{\PSS}$ for rational strongly semi-positive closed symplectic
manifolds is similar to the weakly exact case.
If $(M, \go)$ is not weakly exact, the action functional $\ca_H$
is well-defined only on a suitable infinite cover of
$\cl^\circ$, however.
As a consequence, the set of critical points of $\ca_H$ is
infinite for every $H \in \ch$. 
One therefore needs to define the Floer chain complex over a certain
Novikov ring, and the Floer homology will thus be a module over this ring.
The strong semi-positivity condition excludes bubbling off of
pseudo-holomorphic spheres in the compactifications of the
moduli spaces relevant to the definition of Floer homology,
and it is also used in the construction of the PSS-isomorphism.
The selector $\gs_{\PSS}$ can then be defined as before. 
Since $(M,\go)$ is rational, the spectrum $\Sigma^\circ (H)$ is
a closed and nowhere dense subset of $\RR$ for all $H \in \ch$,
and so the axioms (AS1), (AS3), (AS4) and (AS5) for a weak action
selector can be verified as in the weakly exact case. 
The verification of (AS2), however, makes use of the Poincar\'e duality for
Floer homology.
It is available only if $M$ is closed, and thus no weak action
selector has yet been constructed for convex strongly semi-positive
symplectic manifolds.

\section{Almost existence via $c_{\HZ}$ and $c_{\HZ}^\circ$}  \label{app:1}
\ni
In this appendix we explain why the finiteness of the symplectic
capacities $c_{\HZ}$ and $c_{\HZ}^\circ$ leads to almost
existence results for periodic orbits.
We shall only look at $c_{\HZ}^\circ$, the argument for
$c_{\HZ}$ being the same.
Let $(M, \go)$ be a symplectic manifold. For each $\eps \in
[0,1]$ and each $A \subset M$ we consider the set $\cf^\eps(A)$
of functions in $\ch (A)$ satisfying
\begin{itemize}
\item[(P1)\,\,] 
  $H \ge 0$,
\s
\item[(P2)\,\,] 
  $H |_U = \max H$ for some open non-empty set $U \subset A$,
\s
\item[(P$3^\eps$)] 
  the critical values of $H$ lie in $[0, \eps \max H] \cup \max H$.
\end{itemize}
We denote by $\cf_{\HZ}^{\circ,\eps}(A,M)$ the set of those $H \in
\cf^\eps(A)$ for which the flow $\gf_H^t$ has no non-constant,
contractible in $M$, $T$-periodic orbits
with period $T \le 1$, and set
\[
C_{\HZ}^{\circ,\eps} (A,M) \,=\, 
\sup \left\{ \max H \mid H \in
\cf_{\HZ}^{\circ,\eps} (A,M) \right\} .
\]
Notice that $C_{\HZ}^{\circ,1} (A,M) = C_{\HZ}^{\circ} (A,M)$
and $C_{\HZ}^{\circ,0} (A,M) = c_{\HZ}^{\circ} (A,M)$. Moreover,
we claim that for every $\eps \in (0,1)$,
\begin{equation}  \label{e:1eps}
(1-\eps) C_{\HZ}^{\circ,\eps} \left( A,M \right) \,\le\, c_{\HZ}^\circ
\left( A,M \right) \,\le\, C_{\HZ}^{\circ,\eps} \left( A,M \right).
\end{equation}
The second inequality follows from definitions.
To prove the first inequality,
we need to show that for each $H \in \cf_{\HZ}^{\circ, \eps} (A,M)$ and
each $\delta >0$ there exists $K \in
\cs_{\HZ}^\circ (A,M)$ with $\max K = (1-\eps) \max H -\delta$.
The idea of the argument is to take as $K$ only ``the upper part of $H$''. 
To be more precise, fix $H \in \cf_{\HZ}^{\circ, \eps} (A,M)$ and
$\delta \in \;]0,(1-\eps) \max H [$. 
Following \cite[Section~6]{GG} and \cite{Schl}, we choose a surjective 
smooth map 
$f \colon [0,\max H] \ra [0,(1-\eps)\max H -\delta]$ such that
\[
f(t) = 0 \;\text{ if }\, t \in [0,\eps \max H]
\quad \text{ and } \quad
0 \le f'(t) \le 1 \;\text{ for all }\, t \in [0,\max H] .
\]
Set $K = f \circ H$. Since $|f'|\leq 1$ and since the flow of $H$ has no non-trivial
contractible periodic orbits with period $T\leq 1$, the same holds for
the flow of $K$. Furthermore,
since $H \in \cf_{\HZ}^{\circ, \eps} (A,M)$, we have 
$K \in \cs_{\HZ}^\circ (A,M)$, and by construction $\max K = (1-\eps) \max H -\delta$,
which completes the proof of \eqref{e:1eps}.
\begin{corollary}  
Assume that $A$ is a subset of a symplectic manifold $(M, \go)$
such that $c_{\HZ}^\circ (A,M) < \infty$.
Then almost every regular energy level of a proper function on
$A$ carries a periodic orbit which is contractible in $M$.
\end{corollary}  
\proof
The corollary is proved for $C_{\HZ}^\circ$ in \cite{HZ, MaS},
and the argument given there applies to $C_{\HZ}^{\circ,\eps}$ for every
$\eps \in (0,1]$. The corollary thus follows from the first
inequality in \eqref{e:1eps}.
\proofend

\section{Proofs of Propositions~\ref{p:spec} and \ref{p:restricted}}  \label{app:2}
\ni
{\bf Proof of Proposition~\ref{p:spec}.}
(i) 
If $S$ is simply connected, we choose a contact form $\gl$ on
$S$;
for $x \in \cp^\circ (S)$ we can choose $\bar{x} \in \cd (x)$
contained in $S$, so that $\int_{\bar x} \go = \int_x \gl$.
If $S$ is of restricted contact type, we choose a globally
defined contact form $\gl$; for $x \in \cp^\circ  (S)$ we find again
that $\int_{\bar x} \go = \int_x \gl$.
The Reeb vector field $R$ on $S$ associated with $\gl$ is
defined by
\[
\iota_R \go = 0 
\quad \text{ and } \quad
\gl (R) =1 .
\]
Choose a Riemannian metric on $S$ such that the length of $R$ is
$1$. Parametrizing $x$ such that $\dot x (t) =R$, we see that
$\int_x \gl$ is the length of $x$ and is
hence positive.
If $x_j$ is a sequence in $\cp^\circ (S)$ with $\int_{x_j} \gl \ra
a$, then the lengths $\int_{x_j} \gl$ are bounded, and
hence a subsequence of $x_j$ converges to a closed characteristic
$x \in \cp^\circ (S)$ of length $a$, see \cite[p.\ 109]{HZ}. This
proves that $\Sigma^\circ (S)$ is closed. 
The proof of (ii) is similar.
\proofend

\ni
{\bf Proof of Proposition~\ref{p:restricted}.}
We prove the first statement in Proposition~\ref{p:restricted};
the second statement can be proved in a similar way.
The inequality $0< \ga_1^\circ (S)$ follows from
Proposition~\ref{p:spec}\,(i). We will show 
\begin{equation}  \label{est:ac}
\ga_1^\circ (S) \,\le c_{\HZ}^\circ \left( U,M \right) 
\end{equation} 
for the capacity $C_{\HZ}^\circ (U,M)$. 
It will then be clear from the proof that \eqref{est:ac} holds for
all capacities $C_{\HZ}^{\circ,\eps}(U,M)$, $\eps \in (0,1]$,
introduced in Appendix~\ref{app:1},
and hence also for $c_{\HZ}^\circ(U,M)$ in view of \eqref{e:1eps}.

\s
\ni
{\bf Step 1.}  
{\it A convenient thickening of $S$.}
Let $X$ be a Liouville vector field on $M$ transverse to $S =
\pp U$.
Then $X$ points outward.
Let $\gl = \iota_X \go |_S$ be the associated contact form on
$S$ and let $R$ be the Reeb vector field of $\gl$ on $S$.
Using the flow $\gf_X^t$ of $X$, which exists in a neighborhood of
$\overline{U}$ for small $t$, we see that
a neighborhood of $S$ is symplectomorphic to the thickening
$(-\eps, \eps) \times S$ with coordinates $(t,x)$ and symplectic
form 
\begin{equation}  \label{id:otx}
\go (t,x) \,=\, d \left( (1+t) \gl (x) \right) .
\end{equation}
Then the Liouville vector field is $X(t,x) = (1+t)
\frac{\pp}{\pp t}$.
Set $S_t = \{ t \} \times S$.
Using \eqref{def:H} and \eqref{id:otx}, we see that 
the Hamiltonian vector field $X_H$ of $H \colon (-\eps, \eps)
\times S$, $H(t,x) = t$,
points along $S_t$ and equals the {\it translate}\, of $R$ on
$S_t$. Thus,  the periodic orbits of the flow of $H$ on $S_t$
are exactly the periodic orbits of the Reeb flow on $S$.
Set $U_t = U \cup (-t,t) \times S$.
Then $U_t = \gf_X^{\ln (1+t)} (U)$, and hence the conformality of
the symplectic capacity $C_{\HZ}^\circ$ implies that
\begin{equation}  \label{id:1+t}
C_{\HZ}^\circ \left( U_t, M \right) \,=\, (1+t)\, C_{\HZ}^\circ
\left( U,M \right). 
\end{equation} 

\ni
{\bf Step 2.}
Abbreviate $C = C_{\HZ}^\circ \left( U,M\right)$.
Arguing by contradiction, we assume that $\ga_1^\circ (S) > C$.
Let $\Sigma^\circ (R)$ be the set of periods of those periodic
orbits of the Reeb field $R$ on $S$ which are contractible in
$M$.
Since $\Sigma^\circ (S) = \Sigma^\circ (R)$, the number $\ga_1^\circ (S)$
is the infimum of $\Sigma^\circ (R)$, and so we find $L>1$ such
that 
\begin{equation}  \label{id:empty}
[0,LC] \cap \Sigma^\circ (R) \,=\, \emptyset .
\end{equation}
Fix a small $\gd >0$. Using \eqref{id:empty}, we find $K \in
\cf_{\HZ}^\circ (U,M)$ such that
\[
\max K \,=\, (1-\eps-\gd)\, C.
\]
Let $F$ be a smooth function on the shell $(-\eps,0) \times S$
such that $F(x,t) = f(t)$, where $f \colon (-\eps, 0) \ra \RR$
is a monotone decreasing function such that
\begin{eqnarray*}  
 \begin{array}{rcll} 
    f(t) &=& l \eps C & \text{ for $t$ near $-\eps$}, \\[0.4em]
    f(t) &=& 0        & \text{ for $t$ near $0$}, \\[0.4em]
    | f'(t) | &<& LC  & \text{ for all $t$},  
 \end{array}
\end{eqnarray*}
where $l>1$.
We extend $F$ to the ambient manifold $M$ in the obvious way: $F
\equiv l \eps C$ inside the shell and $F \equiv 0$ outside the shell.
By the construction of the thickening $(-\eps,\eps) \times S$, by
\eqref{id:empty} and by the choice of $f$, the function $F$
belongs to $\cf_{\HZ}^\circ (U,M)$.
Since $F \equiv l\eps C$ on the support of $K$, the Hamiltonian
$H = K+F$ also belongs to $\cf_{\HZ}^\circ (U,M)$.
However,
\[
\max H \,=\, (1-\eps-\gd + l \eps)\, C.
\]
Since $l>1$, we can choose $\gd >0$ so small that $\max H >C$.
This contradicts $H \in \cf_{\HZ}^\circ (U,M)$.
\proofend
\begin{remark}
{\rm
In the above proof, the assumption that $S$ is of {\it
restricted}\, contact type was used only to obtain
identity~\eqref{id:1+t}.
If $S$ is a hypersurface of contact type bounding $U$ such that
the functions $t \mapsto c_{\HZ}^\circ (U_t,M)$ and $t \mapsto
c_{\HZ} (U_t)$ are $1$-Lipschitz at $0$, then
Proposition~\ref{p:restricted} still holds for $S$.
}
\end{remark}

\enddocument